\documentclass[11pt]{amsart}
\usepackage{amsmath,amsfonts,amscd,amssymb}
\newtheorem{theorem}{Theorem}
\newtheorem{lemma}{Lemma}
\newtheorem{corollary}{Corollary}
\theoremstyle{remark}\newtheorem{remark}{Remark}
\newcommand{\field}[1]{\ensuremath{\mathbb{#1}}}
\newcommand{\CC}{\field{C}}
\newcommand{\HH}{\field{H}}

\newcommand{\RR}{\field{R}}
\newcommand{\del}{\partial}

\newcommand{\vphi}{\varphi}
\DeclareMathOperator{\PSL}{PSL} \DeclareMathOperator{\PSU}{PSU}
\DeclareMathOperator{\Imm}{Im}
\begin{document}
\title[Hyperbolic 2-spheres, accessory parameters and
K\"{a}hler metrics]
{Hyperbolic $2$-spheres with conical singularities, accessory parameters\\
and K\"{a}hler metrics on ${\mathcal{M}}_{0,n}$}
\author{Leon Takhtajan and Peter Zograf}
\address{Department of
Mathematics SUNY at Stony Brook, Stony Brook NY 11794-3651 USA.}
\email{leontak@math.sunysb.edu}
\address{Steklov Mathematical Institute, St.Petersburg, 191011 Russia}
\email{zograf@pdmi.ras.ru}

\thanks{Research of the first author was partially supported by the NSF
grant DMS-9802574.}

\begin{abstract}
We show that the real-valued function $S_\alpha$ on the moduli space
${\mathcal{M}}_{0,n}$ of pointed rational curves, defined as the critical value
of the Liouville action functional on a hyperbolic $2$-sphere with 
$n\geq 3$ conical singularities of arbitrary orders $\alpha=\{\alpha_1,\dots,
\alpha_n\}$, generates accessory parameters of the associated 
Fuchsian differential equation as their common antiderivative. 
We introduce a family of K\"ahler metrics on ${\mathcal{M}}_{0,n}$ 
parameterized by the set of orders $\alpha$, explicitly relate accessory 
parameters to these metrics, and prove that the functions $S_\alpha$ are 
their K\"{a}hler potentials.
\end{abstract}

\maketitle

\section{Introduction}  The existence and uniqueness of a hyperbolic metric
(a conformal metric of constant negative curvature $-1$) with prescribed
singularities at a finite number of points on a Riemann surface is a
classical problem that is closely related (and in special cases is
equivalent) to the famous Uniformization Problem of Klein and Poincar\'e.
Actually, in 1898 Poincar\'e \cite{Poin} solved this problem for the simplest
case of {\em parabolic} singularities. Below we formulate his result for the
particular case of the standard 2-sphere
realized as the Riemann sphere $\widehat{\CC}=\CC\cup\{\infty\}$.
Consider the punctured surface $X=\widehat{\CC}\setminus\{z_1,\dots,z_n\}$
with $n\geq 3$ (by
applying an appropriate M\"obius transformation we can always assume
that $z_{n-2}=0, z_{n-1}=1, z_n=\infty$). Then the Liouville  equation
\begin{equation*}
\vphi_{z\bar z}=\frac{1}{2}\,e^{\vphi},
\end{equation*}
(where subscripts stand for the corresponding partial derivatives)
has a unique (real-valued) solution $\varphi$ on $X$ with the following
asymptotics:
\begin{equation*}
\varphi(z)=
\begin{cases}
-2\log|z-z_i| -2\log|\log|z-z_i|| + O(1)& \text{as}~z\rightarrow z_i,\; i\neq n, \\
-2\log|z| - 2\log\log|z| + O(1)& \text{as}~z\rightarrow\infty
\end{cases}
\end{equation*}
(such a singularity is called parabolic).
Geometrically, the Liouville equation mean that the conformal metric
$ds^2=e^\vphi |dz|^2$
on $X$ has constant negative curvature $-1$ (that is, hyperbolic), and
the above asymptotics of $\vphi$ guarantee that $ds^2$ is complete,
and the area of $X$ is $2\pi(n-2)$.

Poincar\'{e} used this result to prove the uniformization theorem, i.e. to show that
there exists a complex-analytic covering of the Riemann surface $X$ by the upper half-plane
$\HH=\{z\in\CC\,|\Imm z>0\}$. He introduced the quantity
\begin{equation*}
T_\vphi = \vphi_{zz} - \frac{1}{2}\, \vphi_z^2
\end{equation*}
and showed that when $\vphi$ satisfies the Liouville equation with parabolic
singularities, then $T_\vphi$ is a meromorphic function on $\widehat{\CC}$ of
the form
\begin{equation*}
T_\vphi(z)=\sum_{i=1}^{n-1}\left(\frac{1}{2(z-z_i)^2}\, + \frac{c_i}{z-z_i}\right),
\end{equation*}
with the asymptotics
\begin{equation*}
T_\vphi(z)=\frac{1}{2z^2}\, + \frac{c_n}{z^3}\, +
O\left(\frac{1}{z^4}\right)~\text{as}~z\rightarrow\infty.
\end{equation*}
The coefficients $c_i$ are the famous accessory parameters. They satisfy
three obvious linear relations imposed  by the asymptotic behaviour of $T_\vphi$
at $\infty$. The coefficients $c_1,\dots,c_n$ are uniquely characterized
by the fact that the monodromy group of the Fuchsian differential
equation
\begin{equation*}
\frac{d^2 u}{dz^2} + \frac{1}{2}\, T_\vphi(z) u =0
\end{equation*}
is conjugate in $\PSL(2,\CC)$ to the group of deck transformations of a
covering $\HH\rightarrow X$.

These ideas of Poincar\'{e} got in the spotlight once again about 20 years ago
due to Polyakov's path integral formulation of the bosonic string \cite{Pol} and the
conformal field theory of Belavin-Polyakov-Zamolodchikov \cite{bpz}.
Briefly,  in the quantum Liouville theory the quantity $T_\vphi$ plays the role of the
$(2,0)$-component of the stress-energy tensor  that
satisfies conformal Ward identities reflecting conformal symmetry of the theory.
At the semi-classical level, as it was first observed by Polyakov, the Ward identity
establishes (at the physical level of rigor) a non-trivial relation between the
accessory parameters and the critical value of the Liouville action functional
(see~\cite{Tak1} for details).

In our paper~\cite{Z-T2}, we rigorously proved Polyakov's conjecture using
the Ahlfors-Bers theory of quasiconformal mappings,
and derived simple explicit formulas connecting
the Liouville equation with accessory parameters and the Weil-Petersson metric
on Teichm\"{u}ller space. More specifically, let
\begin{equation*}
\mathcal{Z}_n=\left\{(z_1,\dots,z_{n-3})\in\CC^{n-3}\,|\, z_i\neq 0,1~\text{and}~
z_i\neq z_k~\text{for}~i\neq k\right\}
\end{equation*}
be the configuration space of singular points ($\mathcal{Z}_n$ is isomorphic
to the moduli space $\mathcal{M}_{0,n}$ of $n$-pointed rational curves over $\CC$).
Then there exists a smooth function
$S:\mathcal{Z}_n\rightarrow\RR$ (critical value of the Liouville action functional;
cf. Section 3) such that
\begin{equation} \tag{I}
c_i=-\frac{1}{2\pi}\frac{\del S}{\del z_i},\qquad i=1,\dots,n-3
\end{equation}
and
\begin{equation} \tag{II}
\frac{\partial c_i}{\partial \bar z_k} = \frac{1}{2\pi}\left\langle
\frac{\partial}{\partial z_i}, \frac{\partial}{\partial
z_k}\right\rangle_{WP},\qquad i,k=1,\dots,n-3,
\end{equation}
where $\langle~,~\rangle_{WP}$ denotes the Weil-Petersson metric on
$\mathcal{Z}_n\cong \mathcal{M}_{0,n}$.\footnote{In~\cite{Z-T3}
we formulated and proved analogs of (I)-(II)
for compact Riemann surfaces of arbitrary genus.}
An immediate corollary of (I) and (II) is that the critical value $S$ of the Liouville
action is a potential for the Weil-Petersson metric.\footnote{These results
were used by the second author in the study of the asymptotic behaviour
of accessory parameters for degenerating Riemann surfaces~\cite{Z}.}

Though our methods generalize verbatim to hyperbolic 2-spheres with
{\em elliptic} singularities of finite order (in which case there exists
a ramified covering $\HH\rightarrow \widehat\CC$
branched over singular points $z_1,\dots,z_n$), they no longer work for
conical singularities of general type (see Section 2 for precise definitions).
However, exact analogs of formulas (I) and (II) hold in this general case
as well, provided the orders $\{\alpha_1,\dots,\alpha_n\}$ of singularities
$z_1,\dots,z_n$ satisfy some rather mild natural conditions. Physical
consideration based on semi-classical limits of conformal Ward
identities also suggests the validity of these formulas in a
general situation.

The objective of this paper is to give straightforward proofs of (I)-(II)
in the case of hyperbolic 2-spheres with conical singularities of general
type. Section 2 contains
the definitions and background material about the classical Liouville
equation, including detailed asymptotics of its solution.
In Section 3 we present the action functional
for the Liouville equation, introduced in~\cite{Tak2}, and prove an analogue of
formula (I) --- Theorem 1.\footnote{Recent physicists' paper \cite{CMS} gives
a different, computationally more involved
proof of Theorem 1.} In Section 4 we prove an analogue of formula (II)
that relates accessory parameters to certain K\"ahler metrics on the moduli
space $\mathcal M_{0,n}$ similar to the Weil-Petersson metric --- Theorem 2.
It is worth to notice that the proofs are considerably simpler than those in
\cite{Z-T2} and do not use Teichm\"uller theory.


\section{Background material}

Consider the Riemann sphere $\widehat\CC = \CC\cup\infty$ with
$n \geq 3$ distinct marked points $z_1,\dots,z_n$. As in the Introduction,
we normalize the last three points to be 0, 1 and $\infty$
respectively, so in the sequel we will always assume that $z_{n-2}=0, z_{n-1}=1,
z_n=\infty$. Let $\alpha=\{\alpha_1,\dots,\alpha_n\}$ be a set of
real numbers such that $\alpha_i < 1,\; i=1,\dots,n$, and
\begin{equation} \label{Cond}
\sum_{i=1}^n \alpha_i > 2.
\end{equation}
According to the classical result of Picard \cite{Pic1,Pic2} (see also \cite{Lich}
and, for a modern proof, \cite{Troy})\footnote{It is very instructive to
compare the approaches of \cite{Pic1,Pic2}, \cite{Lich} and \cite{Troy}.}
there exists a unique conformal metric of constant curvature $-1$, or the
\emph{hyperbolic metric}, on $\widehat\CC$ with conical singularities
of order $\alpha_i$ at $z_i,\;i=1,\dots,n$.
Precisely, it means that such a metric has the form $ds^2=e^\vphi |dz|^2$,
where $\vphi$ is a smooth function on $X=\CC\setminus\{z_1,\dots,z_{n-1}\}$
satisfying the Liouville equation
\begin{equation} \label{Leq}
\vphi_{z\bar z}=\frac{1}{2}\,e^{\vphi}
\end{equation}
and having the following asymptotics near the singular points:
\begin{equation} \label{as}
\vphi(z)=
\begin{cases}
-2\alpha_i\log|z-z_i| + O(1)& \text{as}~z\rightarrow z_i,\; i \neq n,
\\
-2(2 - \alpha_n)\log|z| + O(1)& \text{as}~ z\rightarrow\infty.
\end{cases}
\end{equation}
The point $z_i$ is then called a \emph{conical singularity} of order $\alpha_i$,
or of angle $\theta_i = 2\pi(1 - \alpha_i)$
(we have $\theta_i>2\pi$ when $\alpha_i<0$).
\begin{remark} If $\alpha_i=1$, then $z_i$ is a parabolic point, or \emph{cusp}
(conical singularity of zero angle),
and the asymptotics \eqref{as} should be replaced by
the one mentioned in the Introduction.
\end{remark}

The configuration space ${\mathcal Z}_n$ of singular points
is an open subset in $\CC^{n-3}$:
\begin{equation*}
{\mathcal Z}_n=\left\{(z_1,\dots,z_{n-3})\in\CC^{n-3} |\, z_i\neq 0,1~\text{and}~
z_i\neq z_k~\text{for}~i\neq k\right\},
\end{equation*}
and is isomorphic to
the moduli space ${\mathcal M}_{0,n}$ of
$n$-pointed rational curves over $\CC$.
For any fixed set of orders $\alpha$ the solution $\vphi$ to the Liouville
equation makes sense as a function
of $n-2$ complex variables $z,z_1,\dots,z_{n-3}$, defined on the space
$$\mathcal{Z}_{n+1}=\left\{(z,z_1,\dots,z_{n-3})\in\CC^{n-2}|\,z, z_i\neq 0,1;\,
z\neq z_i;\,z_i\neq z_k~\text{for}~i\neq k\right\}.$$
The space $\mathcal{Z}_{n+1}$ is fibered over $\mathcal{Z}_n$ by ``forgetting''
the first coordinate $z$:
the fiber over a point $(z_1,\dots, z_{n-3})\in\mathcal{Z}_n$ is the surface
$\CC\setminus\{z_1,\dots, z_{n-3},0,1\}$.
It follows from the
results of~\cite{Pic2,Lich,Troy}, that $\vphi$ is a real-analytic function on
$\mathcal{Z}_{n+1}$.

The $(2,0)$-component of the stress-energy tensor in the Liouville theory
is given by the expression
\begin{equation} \label{set}
T_\vphi = \vphi_{zz} - \frac{1}{2}\, \vphi_z^2.
\end{equation}
The following result is classical.
\begin{lemma} \label{1}
Let $\vphi$ be the solution to the Liouville equation with conical singularities
\eqref{as}. Then $T_\vphi$ is a meromorphic function on $\widehat\CC$ with
second order poles at $z_1,\dots,z_n$. Explicitly,
\begin{equation}\label{ap}
T_\vphi(z)=\sum_{i=1}^{n-1} \left(\frac{h_i}{2(z-z_i)^2}\, + \frac{c_i}{z-z_i}\right),
\end{equation}
and
\begin{equation}\label{infty}
T_\vphi(z)=\frac{h_n}{2z^2} + \frac{c_n}{z^3} +
O\left(\frac{1}{z^4}\right)~\text{as}~z\rightarrow\infty,
\end{equation}
where $h_i=\alpha_i(2-\alpha_i),\,i=1,\dots,n$.\footnote{The coefficients $h_i$ are
conformal weights in quantum Liouville theory~\cite{Tak2}.}
\end{lemma}

Complex numbers $c_i$
are called \emph{accessory parameters}. They are uniquely determined by the singular
points $z_1,\dots,z_n$ and the set of orders $\alpha$. Formula \eqref{infty}
imposes three linear equations on the parameters $c_1,\dots,c_n$:
\begin{equation*}
\sum_{i=1}^{n-1}c_i=0,\qquad \sum_{i=1}^{n-1}(h_i+2c_i z_i)=h_n,\qquad
\sum_{i=1}^{n-1}(h_iz_i+c_i z^2_{i})=c_n,
\end{equation*}
so that $c_{n-2}, c_{n-1}$ and $c_n$ are explicit linear
combinations of $c_1,\dots,c_{n-3}$ with coefficients depending on $z_i$ and
$\alpha_i$. Real analyticity of $\vphi$ implies that the accessory parameters are
also real-analytic functions on $\mathcal{Z}_n$.

To study the behaviour of $\vphi$ near the singular points more
thoroughly, consider the Fuchsian differential equation
\begin{equation} \label{Fuchs}
\frac{d^2 u}{dz^2} + \frac{1}{2}\, T_\vphi(z) u=0;
\end{equation}
with regular singular points at $z_1,\dots,z_n$. A classical result (see,
e.g.~\cite{Poin}), which follows from the fact that $e^{-\vphi/2}$ is
a solution to \eqref{Fuchs},
asserts that the monodromy group $\Gamma$ of the differential equation
\eqref{Fuchs} is, up
to a conjugation in $\PSL(2,\CC)$, a subgroup of $\PSL(2,\RR)$ (see, e.g.,
\cite {Ford}, \cite{Bir}, or \cite{Kuga}).\footnote{Among many available
references, \cite {Ford} is
classical, \cite{Bir} gives a detailed exposition
of Fuchsian differential equations, and \cite{Kuga} is a  modern introduction
to the subject.} Such a group $\Gamma$
is discrete in $\PSL(2,\RR)$ if and only if $\alpha_i=1-1/l_i$ for all $i=1,\dots,n$,
where $l_i$ is a positive integer or $\infty$.

In case of general conical singularities the monodromy group $\Gamma$
is no longer discrete in $\PSL(2,\RR)$. It is
generated by local monodromies around regular singular points $z_i$,
which, in general, are elliptic
elements $\gamma_i$ of infinite order.  If we denote the  fixed points
of $\gamma_i$ by $w_i, \bar{w}_i$, then
\begin{equation*}
\frac{\gamma_i (z) -w_i}{\gamma_i (z) -\bar{w}_i}=\lambda_i\,\frac{z -w_i}{z
-\bar{w}_i}, \quad i=1,\dots,n,
\end{equation*}
where $\lambda_i=e^{2\pi\sqrt{-1}(1-\alpha_i)}$ is called the multiplier of $\gamma_i$.
\begin{remark}
It is an outstanding problem to find a geometric meaning of the monodromy group
$\Gamma$ in the case of general conical singularities, thus providing another
interpretation for the accessory parameters.
Perhaps, this problem should be considered in the context of
A.~Connes~\cite{Con} non-commutative differential geometry where such group actions
naturally appear.
\end{remark}

Let $w = u_1/u_2$ be the ratio of two linearly independent solutions $u_1, u_2$ of
the differential equation \eqref{Fuchs}. It is a multi-valued meromorphic function on
$\widehat\CC$ with ramification points $z_1,\dots,z_n$, and it is single-valued
on the universal cover of $X=\widehat\CC\setminus\{z_1,\dots,z_n\}$.
It is a classical result of Schwarz that
\begin{equation} \label{Schwarz}
T_\vphi = {\mathcal S}(w)
\end{equation}
on $X$, where ${\mathcal S}(w)$ denotes the Schwarzian derivative of $w$:
\begin{equation*}
{\mathcal S}(w) = \frac{w'''}{w'} - \frac{3}{2}\left(\frac{w''}{w'}\right)^2.
\end{equation*}

Next, normalize $u_1, u_2$ in such a way that the monodromy group $\Gamma$
of \eqref{Fuchs} is a subgroup of
$\PSU(1,1)$. The multi-valued function $w$ admits the following expansion in the
neighborhood of each singular point $z_i$:
\begin{equation} \label{expand}
\sigma_i(w(z))  = \zeta_{i}^{1-\alpha_i}\sum_{k=0}^\infty a^{(k)}_i \zeta_{i}^k
\quad\text{as}~\zeta_i\rightarrow 0,\;i=1,\dots,n.
\end{equation}
Here $\zeta_i$ is a local uniformizer:
$\zeta_i=z-z_i~\text{for}~i=1,\dots,n-1,~\text{and}~\zeta_n=1/z$, and
$\sigma_i\in\PSU(1,1)$ diagonalizes local monodromy $\gamma_i$
around $z_i$, $i=1,\dots,n$.
Moreover, the coefficients $a^{(k)}_i$ are (locally) real-analytic on $\mathcal{Z}_n$,
as it follows from the analytic dependence on parameters of solutions to ordinary
differential equations.

\begin{lemma} \label{2}
The solution $\vphi$ to the Liouville equation \eqref{Leq}
with conical singularities \eqref{as} is given
by the formula
\begin{equation*}
e^\vphi=\frac{4|w'|^2}{(1-|w|^2)^2}\,,
\end{equation*}
where $w=u_1/u_2$, and $u_1, u_2$ are two linearly independent solutions
of the Fuchsian differential equation \eqref{Fuchs} with monodromy in $\PSU(1,1)$.
\end{lemma}
\begin{proof}
Since the monodromy is in $\PSU(1,1)$, the function
$\log\left(\frac{4|w'|^2}{(1-|w|^2)^2}\right)$ is real and single-valued on $X$.
Moreover, it is easy to check that this function satisfies the
Liouville equation, and by \eqref{expand} it has the same asymptotics \eqref{as}
as $\vphi$. Therefore, it must be equal to $\vphi$.
\end{proof}

\begin{remark} When $\alpha_i=1,\,i=1,\dots,n$, it is more convenient to
normalize solutions $u_1, u_2$ so that $\Gamma\subset\PSL(2,\RR)$ (see~\cite{Z-T2}).
\end{remark}

From the equality \eqref{Schwarz} and expansions \eqref{expand} we readily get the
following formula for the accessory parameters (cf.~Lemma 1 in~\cite{Z-T2}).
\begin{lemma} \label{3}
\begin{equation*}
c_i=\frac{h_i}{1-\alpha_i}\, \cdot \frac{\;a^{(1)}_i}{\;a^{(0)}_i}, \quad i=1,\dots,n,
\end{equation*}
where $h_i=\alpha_i(2-\alpha_i)$.
\end{lemma}

Finally, we summarize all the necessary facts about the asymptotic behaviour
of $\vphi$ and its derivatives in the next statement (cf.~Lemma
2 in~\cite{Z-T2}).

\begin{lemma} \label{4} The solution $\vphi$ to the Liouville equation
\eqref{Leq} with conical
singularities \eqref{as} has the following asymptotic expansions near the singular
points $z=z_i$, uniform in a neighborhood of $(z_1,\dots,z_{n-3})$ in
$\mathcal{Z}_n$:
\begin{itemize}

\item[(i)]
\begin{equation*}
\vphi_z(z)=
\begin{cases}
-\,\dfrac{\alpha_i}{\zeta_i}\, + \,\dfrac{c_i}{\alpha_i}\, + \,\dfrac{f_i(|\zeta_i|)}
{\zeta_i}\,+\,o(1)\quad\text{as}~z\rightarrow z_i,\;i\neq n,\\ -(2-\alpha_n)\zeta_n \,
- \,\dfrac{c_n}{\alpha_n}\cdot\zeta_{n}^2\, + \, f_n(|\zeta_n|)\zeta_n
\,+\,o\left(|\zeta_n|^{2}\right)\quad\text{as}~z\rightarrow\infty,
\end{cases}
\end{equation*}
where $\zeta_i=z-z_i\; (i\neq n)$ and $\zeta_n=1/z$ are local coordinates near
the singular points, and
$$f_i(t)=O\left(t^{2(1-\alpha_i)}\right)\quad\text{as}~t\rightarrow 0,\;i=1,\dots, n$$

\item[(ii)] For $i=1,\dots,n-3$
\begin{equation*}
\vphi_{zz}(z) = \frac{\alpha_i + g^{(0)}_i(\zeta_i) + \zeta_i g^{(1)}_i (\zeta_i)}
{\zeta^{2}_i} + O(1),
\end{equation*}
where
\begin{equation*}
g^{(0)}_i(t),\;g^{(1)}_i (t) = O\left(t^{2(1-\alpha_i)}\right)
\quad\text{as}~t\rightarrow 0.
\end{equation*}

\item[(iii)] For $i=1,\dots,n-3~\text{and}~k=1,\dots,n$, there exist constants $d_{ik}$
such that
\begin{equation*}
\vphi_{z_i}(z)=
\begin{cases}
-\delta_{ik}\vphi_z(z) + d_{ik} + o(1)&~ \text{as}~z\rightarrow z_k,\;k\neq n, \\
\quad d_{in} + o(1)&~\text{as}~z\rightarrow\infty.
\end{cases}
\end{equation*}

\item[(iv)] If $\alpha_k>0$ for each $k=1,\dots,n$, then for $i=1,\dots,n-3$
\begin{equation*}
-2e^{-\vphi}\vphi_{z_i\bar z}=
\begin{cases}
\delta_{ik} + O\left(|z-z_k|^{\;\min\{1,\;2\alpha_k\}}\right)&~\text{as}~z\rightarrow
z_k,\;k\neq n, \\
O\left(|z|^{\;\max\{1,\;2(1-\alpha_n)\}}\right)&~\text{as}~z\rightarrow \infty.
\end{cases}
\end{equation*}
\end{itemize}
\end{lemma}
\begin{proof}
Parts (i)-(iii) follow from \eqref{expand} and Lemmas \ref{2} and \ref{3}; part
(iv) follows from (i), (iii), the Liouville equation \eqref{Leq}
and asymptotics \eqref{as}.
Uniform estimates for the remainder terms follow from the real analyticity
of the coefficients $a^{(k)}_i$ as functions of $z_1,\dots,z_{n-3}$.
One can also prove (i)-(iv) directly from the Liouville
equation and asymptotics \eqref{as} by observing that the solution $\vphi$
admits the following expansion in a neighborhood of each $z_i$:
\begin{equation*}
\vphi(z) = -2\alpha_i \log|z-z_i| + \xi^{(0)}(z) + \sum_{k=1}^\infty
|z-z_i|^{\,2k(1-\alpha_i)}\xi^{(k)}(z),\;i=1,\dots,n-1,
\end{equation*}
and a similar expansion at $\infty$, where $\xi^{(k)}(z)$ are
real-analytic as functions on the fibered space $\mathcal Z_{n+1}$
(real-analytic dependence on $z_1,\dots,z_{n-3}$ follows from the analysis
in~\cite{Pic1,Pic2,Lich,Troy}).
\end{proof}

\section{Liouville action and accessory parameters}
For a given set of orders $\alpha=\{\alpha_1,\dots, \alpha_n\}$
the action functional for the Liouville equation (\ref{Leq}) is defined
in~\cite{Tak2} by the formula
\begin{equation}
S_\alpha [\psi]=\lim_{\varepsilon \rightarrow 0}
S_\alpha^\varepsilon[\psi],\label{Lac}
\end{equation}
where
\begin{align}\label{action}
S_\alpha^\varepsilon [\psi] &=
\iint\limits_{X^{\varepsilon}}(|\,\psi_{z}|^{2}+e^{\psi})\, \left|\frac{dz\wedge
d\bar z}{2}\right|\\
& + \frac{\sqrt{-1}}{2}\sum_{i=1}^{n-1} \alpha_i\oint\limits_{C_i^\varepsilon}
\psi\left(\frac{d \bar{z}} {\bar{z}-\bar{z}_i} - \frac{dz}{z-z_i}\right) \nonumber\\
& + \frac{\sqrt{-1}}{2}(2-\alpha_n)
\oint\limits_{C_n^\epsilon}\psi\left(\frac{d\bar{z}}{\bar{z}} - \frac{dz}{z}\right)
\nonumber \\
& -2 \pi\sum_{i=1}^{n-1}\alpha^{2}_{i}\log\varepsilon
-2\pi(2-\alpha_n)^2 \log\varepsilon. \nonumber
\end{align}
Here $X^{\varepsilon}=\CC \setminus
\left(\bigcup^{n-1}_{i=1}\{|z-z_i|<\varepsilon\}\bigcup \{|z|>1/\varepsilon\}\right)$,
and the circles $C_i^\varepsilon = \{ |z-z_i| = \varepsilon
\},\,i=1,\dots,n-1,~\text{and}~C_n^\varepsilon = \{ |z| = 1/\varepsilon \}$ are
oriented as the boundary components of $X^\varepsilon$. The Liouville equation
is the Euler-Lagrange equation for the functional $S_\alpha$,
which is defined on the space of all
conformal metrics $e^\psi\,|dz|^2$ on $\widehat\CC$ with conical singularities
at $z_1,\dots,z_n$ of orders $\alpha_1,\dots,\alpha_n$,
satisfying
\begin{equation} \label{der}
\psi_z(z)=
\begin{cases}
-\dfrac{\alpha_i}{z - z_i}\,(1 + o(1))&~\text{as}~ z\rightarrow z_i,\;i\neq n, \\
-(2-\alpha_n)\,\dfrac{1}{z}\,(1 + o(1))&~\text{as}~z\rightarrow\infty.
\end{cases}
\end{equation}
\begin{remark}
The contour integrals in \eqref{action} ensure that for all $\psi$
satisfying \eqref{as} and \eqref{der} and for all $u\in C^\infty(\widehat\CC,\RR)$
\begin{equation*}
\lim_{t\rightarrow 0}\frac{S[\psi +tu]-S[\psi]}{t}=
\iint\limits_{\CC}(-2\psi_{z\bar z} + e^\psi)\,u\,
\frac{|dz\wedge d\bar z|}{2}.
\end{equation*}
\end{remark}

The Liouville action evaluated on the solution $\vphi$ to the Liouville equation is a real
valued function $S_\alpha[\vphi] = S_\alpha (z_1,\dots,z_{n-3})$ on the configuration
space ${\mathcal Z}_n$ depending on $\alpha_1,\dots,\alpha_n$ as parameters.

\begin{theorem}
For any fixed set of orders $\alpha = \{\alpha_1,\dots,\alpha_n\}$
such that $\alpha_i < 1$ and
$\sum_{i=1}^n \alpha_i > 2$, the function $S_\alpha :{\mathcal
Z}_n \longrightarrow\RR$ is differentiable and
\begin{equation}
c_i=\frac{1}{2\pi}\frac{\partial S_\alpha}{\partial z_i},\qquad i=1,\dots,n-3,
\end{equation}
where $c_i$ are the accessory parameters defined by \eqref{ap}.
\end{theorem}

\begin{proof}
First we show that
\begin{equation}
\lim_{\varepsilon\rightarrow 0}\frac{\partial S_\alpha^\varepsilon}{\partial z_i}=
-2\pi c_i \label{limit}
\end{equation}
pointwise on the configuration space ${\mathcal Z}_n$. We have
\begin{align}\label{ea}
\frac{\partial S_\alpha^\varepsilon}{\partial z_i}&=\frac{\sqrt{-1}}{2}
\left(\iint\limits_{X^\varepsilon}\frac{\partial }{\partial z_i} (|\,\vphi_z|^2 +
e^\vphi)\, dz\wedge d\bar z + \oint\limits_{C_i^\varepsilon}(|\,\vphi_z|^2 +
e^\vphi)\, d\bar z \right)
\\ & + \frac{\sqrt{-1}}{2} \sum_{k=1}^{n-1} \alpha_k \oint\limits_{C_k^\varepsilon}
(\vphi_{z_i} + \delta_{ik}\vphi_z) \left(\frac{d \bar{z}} {\bar{z}-\bar{z}_k} -
\frac{dz}{z-z_k}\right)\nonumber
 \\ & + \frac{\sqrt{-1}}{2}\,(2-\alpha_n)
\oint\limits_{C_n^\epsilon}\vphi_{z_i}\left(\frac{d\bar{z}}{\bar{z}} -
\frac{dz}{z}\right).\nonumber
\end{align}
Using part (i) of Lemma \ref{4}, we see that
$$\frac{\sqrt{-1}}{2} \oint\limits_{C_i^\varepsilon}(|\,\vphi_z|^2 \, d\bar z
\longrightarrow \pi c_i\qquad\text{as}~\varepsilon\rightarrow 0.$$
From the Liouville equation we get
\begin{equation*}
\oint\limits_{C_i^\varepsilon} e^\vphi d\bar z = -\frac{1}{2}
\oint\limits_{C_i^\varepsilon}\vphi_{zz} dz,
\end{equation*}
which tends to $0$ as $\varepsilon\rightarrow 0$ because of part (ii)
of Lemma~\ref{4}. As it follows from part
(iii) of Lemma~\ref{4}, the contour integrals in the second and third lines
of \eqref{ea} tend to
\begin{equation*}
-2\pi c_i -2\pi \sum_{k=1}^{n-1}\alpha_k d_{ik} -2\pi (\alpha_n-2) d_{in}
\end{equation*}
as $\varepsilon\rightarrow 0$. An obvious identity
\begin{equation*}
\frac{\partial}{\partial z_i} |\vphi_z|^2 dz\wedge d\bar z = d\,(\vphi_{z_i}\vphi_{\bar
z}\, d\bar z - \vphi_{z_i}\vphi_z\, dz) - 2\vphi_{z_i}\vphi_{z\bar z}\,dz\wedge d\bar z,
\end{equation*}
combined with the Liouville equation
yields the following simple formula:
\begin{equation} \label{formula}
\frac{\partial}{\partial z_i} (|\vphi_z|^2 + e^\vphi)\,dz\wedge d\bar z =
d\,(\vphi_{z_i}\vphi_{\bar z}\, d\bar z - \vphi_{z_i}\vphi_z\, dz).
\end{equation}
This reduces the area integral in \eqref{ea} to a
sum of contour integrals. These contour integrals are again easy to
evaluate using parts (i)
and (iii) of Lemma~\ref{4}, and all together they tend to
\begin{equation*}
-\pi c_i + 2\pi\sum_{k=1}^{n-1}\alpha_k d_{ik} + 2\pi (\alpha_n-2) d_{in}
\end{equation*}
as $\varepsilon\rightarrow 0$.
Adding all the terms in the right hand side of \eqref{ea}, we get $-2\pi c_i$
in the limit
as $\varepsilon\rightarrow 0$. Finally, we observe that the convergence of
\eqref{limit} is uniform on compact subsets of ${\mathcal Z}_n$
because so are the estimates in Lemma~\ref{4}.
\end{proof}
\begin{remark}
The same method works for $\alpha_i=1,\,i=1,\dots,n$. In this case formula
\eqref{action} for the functional $S^\varepsilon[\vphi]$ contains an
additional regularizing term $4\pi(n-2)\log\log|\varepsilon|$.
By part 2) of Lemma 2 in~\cite{Z-T2},  no contour integrals
contribute to the classical
action $S[\vphi]$. This gives a much simpler proof of Theorem 1 in~\cite{Z-T2}
along the lines of this paper, without using either the uniformization
theorem, or the quasiconformal mappings.
\end{remark}

\section{Accessory parameters and K\"{a}hler metrics on $\mathcal{M}_{0,n}$}
Throughout this section we assume, in addition, that the orders
$\alpha_1,\dots,\alpha_n$ are all positive,\footnote{This
is equivalent to the condition that all conformal
weights $h_i$ are positive.} i.e., $\alpha_i\in (0,1)$ for each
$i=1,\dots,n$, and $\sum_{i=1}^n \alpha_i >2$. To every such
set of orders $\alpha=\{\alpha_1,\dots,\alpha_n\}$ we associate
Hermitian metric on the configuration space $\mathcal{Z}_n\cong
\mathcal{M}_{0,n}$ as follows.

Consider the kernel
\begin{equation} \label{R}
R(\zeta, z)=-\frac{1}{\pi}\left(\frac{1}{\zeta - z} + \frac{z - 1}{\zeta} -
\frac{z}{\zeta - 1}\right),\qquad (\zeta, z)\in\CC\times\CC,
\end{equation}
and put
\begin{equation} \label{q-diff}
Q_i(z)=R(z, z_i),\quad i=1,\dots,n-3.
\end{equation}
Clearly, the functions $Q_i$ are linear independent.
It follows from the positivity of orders $\alpha_i$ and \eqref{as} that
the functions $Q_i$ are square integrable on
$\widehat\CC$ with respect to the measure
$e^{-\vphi}\frac{|dz\wedge d\bar z|}{2}$.
We define the scalar products of the basis of 1-forms on $\mathcal{Z}_n$
over the point $(z_1,\dots,z_{n-3})\in\mathcal{Z}_n$ by the formula
\begin{equation} \label{hermit}
(dz_i, dz_k)_\alpha = \iint\limits_{\CC} Q_i \overline{Q}_k e^{-\vphi}\,\frac{|dz
\wedge d\bar z|}{2},\quad i,k=1,\dots,n-3.
\end{equation}
The scalar products
$\langle\frac{\partial}{\partial z_i}, \frac{\partial}{\partial z_k}
\rangle_{\alpha}$ are given by the elements of the inverse matrix
to  $\{(dz_i, dz_k)_{\alpha}\}_{i,\,k=1}^{n-3}$.
Since the matrix
$\{\langle\frac{\partial}{\partial z_i}, \frac{\partial}{\partial z_k}
\rangle_\alpha\}_{i,\,k=1}^{n-3}$
is non-degenerate and depends real analytically on $z_i$, it gives rise
to a Hermitian metric on ${\mathcal Z}_n$  which we denote by
$\langle\,\cdot\, , \cdot\rangle_{\alpha}$. This metric
is analogous to the celebrated Weil-Petersson metric on
the moduli space ${\mathcal M}_{0,n}$.\footnote{We get the Weil-Petersson
metric if all the orders $\alpha_i$ are equal to 1.}

\begin{remark} In Teichm\"{u}ller theory, when all $\alpha_i=1$, the holomorphic 
cotangent space to ${\mathcal Z}_n$ at the point
$(z_1,\dots,z_{n-3})$ is identified by means of quasiconformal mappings 
with the space of rational functions on
$\widehat{\CC}$ with only simple poles at $z_1,\dots,z_{n-3},0,1,\infty$,
and $dz_i$ then corresponds to $Q_i$ (see, e.g., \cite{Z-T2} and references
therein). Here we use the same identification directly.
\end{remark}

The kernel $R$, roughly speaking, inverts the operator $\partial/\partial\bar z$ on $\CC$.
The precise statement (see, e.g., \cite{Ahl} for details) is essentially a version of
the Pompeiu formula.
\begin{lemma}\label{5}
Let $g$ be a locally integrable function on $\CC$
such that $g(z)=o(|z|)$ as $z\rightarrow\infty$.
Then the equation
\begin{equation*}
f_{\bar z}=g
\end{equation*}
has a unique solution on $\CC$ satisfying $f(0)=f(1)=0$ and $f(z)=
o(|z|^{2})$ as $z\rightarrow\infty$. This solution is explicitly given by the formula
\begin{equation} \label{sol-d-bar}
f(z)=\iint\limits_\CC g(\zeta)R(\zeta,z)\,\frac{|d\zeta\wedge
d\bar\zeta|}{2}.
\end{equation}
\end{lemma}

Let us formulate the main result of this section.
\begin{theorem}
For any set of orders $\alpha=\{\alpha_1,\dots,\alpha_n\}$
such that $\alpha_i\in (0,1)$ for each $i=1,\dots, n$ and
$\sum_{i=1}^n \alpha_i >2$, we have
\begin{equation}\label{dap}
\frac{\partial c_i}{\partial \bar z_k} = \frac{1}{2\pi}\left\langle
\frac{\partial}{\partial z_i}\;, \frac{\partial}{\partial z_k}\right\rangle_\alpha,
\qquad i,k=1,\dots,n-3.
\end{equation}
\end{theorem}

\begin{proof}
As we mentioned in Section 2, the accessory parameters $c_1,\dots,c_{n-3}$
are real-analytic functions on $\mathcal{Z}_n$. Now consider the functions
\begin{equation*}
F^i = -2e^{-\vphi}\vphi_{z_i \bar z},\;i=1,\dots,n-3.
\end{equation*}
According to part (iv) of Lemma \ref{4} we have
\begin{eqnarray}
&& F^i(z_k) = \delta_{ik},\qquad k = 1,\dots,n-1,\label{delta}\\
&& F^i(z) = O(|z|^{\,\max\{1,\,2(1-\alpha_n)\}}), \qquad z\rightarrow\infty.
\nonumber
\end{eqnarray}
Moreover, as it follows from \eqref{set} and \eqref{ap},
\begin{eqnarray*}
F^i_{\bar z}&=&2e^{-\vphi}\vphi_{\bar z}\,\vphi_{z_i \bar z} - 2e^{-\vphi}\vphi_{z_i
\bar z \bar z}\: =\: -2e^{-\vphi}\frac{\partial}{\partial z_i}\left(\vphi_{\bar z \bar
z} - \frac{1}{2}\,\vphi_{\bar z}^2\right)\\ &=&-2e^{-\vphi}\sum_{k=1}^{n-1}\frac{1}{\bar
z - \bar z_k}\,\cdot\frac{\partial \bar c_k} {\partial z_i}\: =\: 2\pi e^{-\vphi}
\sum_{k=1}^{n-3}\frac{\partial \bar c_k}{\partial z_i}\: \overline{Q}_k.
\end{eqnarray*}
Lemma \ref{5} applied to $g=F^i_{\bar z}$ yields
$$F^i(z) = \iint_{\CC}F^i_{\bar \zeta}(\zeta) R\,(\zeta,z)\,
\frac{|d\zeta\wedge d\bar\zeta|}{2}, \qquad i=1,\dots,n-3. $$
Putting $z=z_j$ and using \eqref{delta} we get that
\begin{equation*}
\delta_{ij} = 2\pi \sum_{k=1}^{n-3}\frac{\del\bar c_k}{\partial z_i}\:
(dz_j, dz_k)_\alpha,\qquad i,j = 1,\dots, n-3,
\end{equation*}
which proves the theorem.
\end{proof}

\begin{remark}
The same arguments prove Theorem 2 in~\cite{Z-T2} making
the uniformization theorem and  quasiconformal mappings redundant
also in the case when all $\alpha_i=1$.
\end{remark}

\begin{corollary} For any set $\alpha$ as in Theorem 2
\begin{equation*}
\left\langle\frac{\partial}{\partial z_i}\;,\frac{\partial}{\partial
z_k}\right\rangle_\alpha = - \,\frac{\partial^2 S_\alpha}{\partial z_i\,\partial
\bar{z}_k},\qquad i,k=1,\dots, n-3.
\end{equation*}
That is, the metric $\langle\,\cdot\, , \cdot\rangle_\alpha$ is K\"ahler and
the function $-S_\alpha$ is its real-analytic K\"ahler potential
on $\mathcal{Z}_n$.
\end{corollary}
\begin{proof}
Immediately follows from Theorems 1 and 2.
\end{proof}

\end{document}